\documentclass[final, 3p, times]{elsarticle}%
\usepackage{graphicx}
\usepackage{xcolor}
\usepackage{colortbl}
\usepackage{epstopdf}
\usepackage{amssymb}
\usepackage{amsthm}
\usepackage{amsmath}
\usepackage{color}
\usepackage{dsfont}
\usepackage{mathtools}
\usepackage{extarrows}
\usepackage{hyperref}
\usepackage{bm}
\usepackage{caption}
\usepackage{float}
\usepackage{verbatim}
\usepackage{epsfig}
\usepackage{multirow}
\usepackage{subfig}
\usepackage[section]{placeins}
\usepackage{float}
\usepackage{cleveref}
\usepackage{diagbox}
\usepackage{algorithm}
\usepackage{algpseudocode}
\usepackage{xr}

\newtheorem{rmk}{Remark}[section]

\biboptions{compress}

\journal{Elsevier}

\begin{document}
\begin{frontmatter}
\title{Preserving Conservation Laws in the Time-Evolving Natural Gradient Method via Relaxation and Projection Techniques}

\tnotetext[label1]{The research of Dongling Wang is supported in part by the National Natural Science Foundation of China under grants 12271463 and the Natural Science Foundation of Hunan Province under grants 2026JJ50363.}

\author[XTU]{Zihao Shi}
\ead{202331510139@smail.xtu.edu.cn}
\author[XTU]{Dongling Wang\corref{mycorrespondingauthor}}
\ead{wdymath@xtu.edu.cn}
\cortext[mycorrespondingauthor]{Corresponding author.}

\address[XTU]{Hunan Research Center of the Basic Discipline Fundamental Algorithmic Theory and Novel Computational Methods, National Center for Applied Mathematics in Hunan, School of Mathematics and Computational Science, Xiangtan University, Xiangtan, Hunan 411105, China}

\begin{abstract}
Neural networks have demonstrated significant potential in solving partial differential equations (PDEs). While global approaches such as Physics-Informed Neural Networks (PINNs) offer promising capabilities, they often lack inherent temporal causality, which can limit their accuracy and stability for time-dependent problems. In contrast, local training frameworks that progressively update network parameters over time are naturally suited for evolving PDEs. However, a critical challenge remains: many physical systems possess intrinsic invariants---such as energy or mass---that must be preserved to ensure physically meaningful solutions. This paper addresses this challenge by enhancing the Time-Evolving Natural Gradient (TENG) method, a recently proposed local training framework. We introduce two complementary techniques: (i) a relaxation algorithm that ensures the target solution $u_{\text{target}}$ preserves both quadratic and general nonlinear invariants of the original system, providing a structure-preserving learning target; and (ii) a projection technique that maps the updated network parameters $\theta(t)$ back onto the invariant manifold, ensuring the final neural network solution strictly adheres to the conservation laws. Numerical experiments on the inviscid Burgers equation, Korteweg-de Vries equation, and acoustic wave equation demonstrate that our proposed approach significantly improves conservation properties while maintaining high accuracy.
\end{abstract}


\end{frontmatter}

\section{Introduction}

\subsection{Motivation and Background}

Neural networks have emerged as a prominent mesh-free alternative for solving partial differential equations (PDEs) \cite{r1,du2021evolutional,han2018solving,li2020fourier}. Their principal advantage lies in the use of automatic differentiation to compute differential operators without spatial grids, offering a viable pathway for simulating high-dimensional PDEs. Raissi et al. proposed Physics-Informed Neural Networks (PINNs) \cite{r1}, which embed governing equations directly into the loss function within a static optimization framework. However, a notable limitation arises for time-dependent problems: PINNs lack an explicit mechanism to enforce temporal causality, as they process temporal and spatial coordinates concurrently \cite{jung2024ceens}. This may yield solutions that violate the temporal consistency required by physical processes.

Existing strategies to address this limitation fall into two categories \cite{jung2024ceens, r11, r10}: soft constraints assign time-varying weights to encourage fitting earlier times first \cite{r10}, while hard constraints adopt classical time-stepping concepts \cite{r11}. Both approaches, however, increase computational overhead and may slow the solving process.

Recently, a novel time-local training approach known as the Neural Galerkin method was introduced \cite{bruna2024neural, du2021evolutional}. This method updates network parameters sequentially in time, evolving from static $\theta$ to time-dependent $\theta(t)$. Grounded in the Dirac-Frenkel time-dependent variational principle \cite{lubich2008quantum}, the right-hand side of the governing equation is projected onto the tangent space of the parametric solution manifold, yielding a dynamical system for $\theta(t)$. In this local training framework, learning becomes the time integration of a parameter dynamics system driven by the underlying physics.

Building upon this foundation, Chen et al. \cite{chen2024teng} proposed the Time-Evolving Natural Gradient (TENG) method, which combines the Neural Galerkin scheme with optimization-based time integration. TENG projects solution evolution in physical space onto network parameter space, guiding parameter updates via natural gradients rather than blind search.

Many dynamical systems possess intrinsic invariants---such as energy, mass, or momentum---that reflect fundamental conservation laws. A key insight from classical numerical analysis is that strict preservation of these invariants is paramount for stability and reliability, sometimes even more crucial than precisely satisfying the original equations \cite{HairerWannerLubich2006, boesen2025neural}. For neural-network-based solvers, ensuring that parameterized approximations respect these physical laws is essential for long-term accuracy and physical fidelity.

When applying TENG to time-dependent PDEs with conservative structures, two fundamental difficulties arise:
(i) Target solution lacks conservation. The target solution $u_{\text{target}}$, typically generated by explicit Runge-Kutta methods, often violates quadratic and nonlinear conserved quantities such as energy. This causes the training data itself to deviate from true physical invariants, compromising the learning target from the outset.
(ii) Parameter updates deviate from the invariant manifold. The tangent-space approximation process in TENG does not incorporate conservation constraints into parameter optimization. Consequently, trained network parameters $\theta$ are not confined to the theoretical conservation manifold $\mathcal{M}_I$, fundamentally limiting the model's ability to predict long-term system evolution accurately.

To address these shortcomings, we introduce a dual-strategy enhancement to the TENG framework:

\begin{itemize}
\item \textbf{Relaxation algorithm for target solution:} Following the relaxation approach proposed by Ranocha et al. \cite{ranocha2020relaxation}, we incorporate a relaxation parameter $\gamma$ that scales the time step to ensure $u_{\text{target}}$ preserves desired invariants while maintaining the accuracy order of the original time discretization.

\item \textbf{Projection algorithm for parameter updates:} Drawing on geometric numerical integration principles \cite{HairerWannerLubich2006, schwerdtner2024nonlinear}, we apply a projection technique that maps updated network parameters back onto the conservation manifold $\mathcal{M}_I$, ensuring the final neural network solution strictly adheres to conservation laws.
\end{itemize}

This combined approach---which we denote as RP-TENG (Relaxation-Projection TENG)---provides a structure-preserving learning target and enforces conservation throughout parameter evolution. Numerical experiments on three benchmark problems demonstrate significant improvements in conservation properties while maintaining computational efficiency and accuracy.

The remainder of this paper is structured as follows. Section 2 reviews the mathematical formulation of structure-preserving PDEs and the Neural Galerkin framework. Section 3 introduces the TENG method and presents our relaxation and projection enhancements. Section 4 provides comprehensive numerical experiments on the inviscid Burgers equation, Korteweg-de Vries equation, and acoustic wave equation. Section 5 concludes and discusses future research directions.

\section{Preliminaries}

This section reviews the mathematical modeling of time-dependent PDEs with conservation laws and the foundational concepts of the Neural Galerkin method.

\subsection{Time-Dependent PDEs with Conservation Laws}

Consider a partial differential equation on a bounded domain $\mathcal{X} \subseteq \mathbb{R}^d$. Let $\mathcal{U}$ be a suitable function space such that solutions $u(\cdot, t)$ over the time interval belong to this space. A PDE with physical structure can be written as
\begin{equation}
\begin{aligned}
&\partial_t u(x,t) = f\bigl(x, u(\cdot, t)\bigr), \qquad \forall (x,t) \in \mathcal{X} \times [0, \infty), \\
&u(x,0) = u_0(x), \qquad \forall x \in \mathcal{X},
\end{aligned}
\label{eq:pde}
\end{equation}
where $f$ contains spatial differential operators such as the Laplacian $\Delta$.

Let $I$ be a physical quantity associated with the equation. For a conservative system, $I$ satisfies
\begin{equation}
\frac{d}{dt} I(u(t)) = 0, \quad \forall t \geq 0, \label{eq:conservation}
\end{equation}
for all solutions of the equation. More generally, dissipative systems satisfy $dI/dt \leq 0$. This paper focuses on conservative systems where equality holds identically.
For numerical approximations, it is crucial that discrete solutions preserve such invariants, i.e., $I(u_n) = I(u_0)$ for all $n \geq 1$. To study conservation laws systematically, we consider invariants of the form
\begin{equation}
I_i(u) = \int_{\mathcal{X}} k_i(u(x,t)) \, d\mu, \quad i = 1, \ldots, m, \label{eq:invariant}
\end{equation}
where $k_i$ may be linear, quadratic, or nonlinear functions of $u$, corresponding to linear, quadratic, or nonlinear invariants, respectively.

\subsection{The Neural Galerkin Framework}
\label{sec:neural_galerkin}

We assume that at any time $t$, the solution $u(\cdot,t)$ admits a parameterized representation $\hat{u}(\cdot;\theta(t))$, where the parameters $\theta(t) \in \Theta \subseteq \mathbb{R}^{n_\theta}$ are time-dependent:
\begin{equation}
u(t, x) \approx \hat{u}(\theta(t), x), \qquad \forall (t, x) \in [0, \infty) \times \mathcal{X}. \label{eq:param_ansatz}
\end{equation}
Denote the parameterized solution manifold as $\mathcal{M} = \{\hat{u}(\eta, \cdot) \mid \eta \in \Theta\} \subseteq \mathcal{U}$, with tangent space
\begin{equation}
T_{\hat{u}(\theta(t),\cdot)} \mathcal{M} = \operatorname{span}\{\partial_{\theta_1}\hat{u}, \ldots, \partial_{\theta_{n}}\hat{u}\}.
\end{equation}

Substituting the ansatz \eqref{eq:param_ansatz} into \eqref{eq:pde} yields the residual
\begin{equation}
r_t(x,\theta(t),\dot{\theta}(t)) = \partial_t \hat{u} - f(x,\hat{u}) = \nabla_\theta \hat{u}(\theta(t),x) \cdot \dot{\theta}(t) - f(x,\hat{u}). \label{eq:residual}
\end{equation}
The Dirac-Frenkel variational principle requires the residual to be orthogonal to the tangent space:
\begin{equation}
\langle \partial_{\theta_i} \hat{u}(\theta(t), \cdot), r_t(\theta(t), \dot{\theta}(t)) \rangle_\mu = 0, \quad i = 1, \ldots, n_\theta, \label{eq:dirac_frenkel}
\end{equation}
where $\langle x, y \rangle_\mu = \int_{\mathcal{X}} x \cdot y \, d\mu$,where $\mu$ denotes a probability measure. This yields the dynamical system for network parameters:
\begin{equation}
M(\theta(t)) \dot{\theta}(t) = F(\theta(t)), \label{eq:param_ode}
\end{equation}
with
\begin{align}
M(\theta(t)) = \int_{\mathcal{X}} \nabla_{\theta}\hat{u}^{\top} \nabla_{\theta}\hat{u} \, d\mu, \quad
F(\theta(t)) = \int_{\mathcal{X}} \nabla_{\theta}\hat{u}^{\top} f(\hat{u}) \, d\mu. \label{eq:mass_force}
\end{align}

Equation \eqref{eq:param_ode} resembles a classical ODE system, with numerical stability determined by the invertibility of the mass matrix $M$. Due to inherent symmetries in neural network parameters (e.g., weight scaling symmetry). Since the matrix $M$ is often singular in practical applications, regularization is essential to maintain numerical stability. Complementing these algorithmic efforts, Lubich has established a rigorous theoretical foundation for regularized Neural Galerkin schemes \cite{feischl2024regularized, lubich2025regularized}.

\begin{rmk}
Even if the exact solution $u(x,t)$ can be represented exactly by the network at each time (i.e., there exists $\theta^*(t)$ such that $u(x,t) = \hat{u}(\theta^*(t), x)$), singularity of $M$ may prevent the optimal parameters from forming a continuous path in parameter space. Consequently, smooth temporal evolution via \eqref{eq:param_ode} cannot be guaranteed \cite{finzi2023stable}.
\end{rmk}

Initial parameters $\theta_0$ are obtained by solving
\begin{equation}
\theta_0 = \arg\min_{\eta \in \Theta} \|\hat{u}(\eta,\cdot) - u_0\|^2_{L^2(\mathcal{X})}. \label{eq:initial_fit}
\end{equation}

\subsection{Structure Preservation in the Neural Galerkin Framework}

The Neural Galerkin method can be formulated as an optimization problem based on an inner product $\langle\cdot,\cdot\rangle_{\mathcal{H}}$:
\begin{equation}
\dot{\theta}(t) = \arg\min_{\eta \in \Theta} \| r_t(\theta(t),\eta) \|^2_{\mathcal{H}}. \label{eq:optimization_form}
\end{equation}
Setting the gradient to zero yields \eqref{eq:param_ode} with the appropriate inner product.
This is equivalent to projecting $f$ onto the tangent space $T_{\hat{u}(\theta,\cdot)}\mathcal{M}$ with respect to $\langle\cdot,\cdot\rangle_{\mathcal{H}}$. Denoting this projection operator as $\mathsf{P}_\theta$, we have
\begin{equation}
\partial_t \hat{u} = \mathsf{P}_\theta f. \label{eq:projection_eq}
\end{equation}
If the projection is inexact, i.e., $\mathsf{P}_\theta f = f + \varepsilon(t)$ with $\varepsilon(t) \neq 0$, then for a conserved quantity $I$ satisfying $I'(u) \cdot f(u) = 0$, we obtain
\begin{align}
\frac{d}{dt} I(\hat{u}) &= I'(\hat{u}) \cdot \partial_t \hat{u} = I'(\hat{u}) \cdot (f(\hat{u}) + \varepsilon(t)) = I'(\hat{u}) \cdot \varepsilon(t) \neq 0. \label{eq:conservation_error}
\end{align}
Thus, projection error strictly prevents preservation of invariants. Thus, projection errors strictly preclude invariant conservation, requiring explicit constraints to achieve an invariant-preserving Neural Galerkin framework.

\section{Relaxation-Projection Enhanced TENG Framework}

When constructing parameterized dynamical systems based on the Neural Galerkin method, stability is often constrained by the ill-posedness of the mass matrix $M$. The implicit function theorem suggests that ensuring continuous evolution of the solution manifold with respect to $\theta$ becomes difficult, posing challenges for numerical solutions. To address these issues, Luo et al. proposed TENG \cite{chen2024teng}, a local training framework that we now extend with relaxation and projection techniques.

\subsection{The TENG Framework}

Optimization-Based Time Integration (OBTI) parameterizes the solution in a continuous function space and computes spatial derivatives via automatic differentiation. Discretizing the time derivative via an explicit Euler scheme gives:
\begin{equation}
u_{\text{target}} = \hat{u}_{\theta_n} + \Delta t \, f(\hat{u}_{\theta_n}). \label{eq:target_euler}
\end{equation}
The parameters $\theta_{n+1}$ corresponding to $u_{\text{target}}$ are obtained by solving
\begin{equation}
\theta_{n+1} = \arg\min_{\theta \in \Theta} L(\hat{u}(\theta,\cdot), u_{\text{target}}), \label{eq:obti}
\end{equation}
where typical loss functions include $L^2$ distance, $L^1$ distance, or KL divergence.

While OBTI provides the best approximation in parameter space, the optimization problem is typically non-convex and challenging to solve directly. The Neural Galerkin scheme, following the Time-Dependent Variational Principle (TDVP), formulates parameter evolution as a least-squares problem:
\begin{equation}
\dot{\theta}(t) = \arg\min_{\dot{\theta}(t)} \| f(\hat{u}_{\theta}) - \nabla_\theta \hat{u}_{\theta} \dot{\theta}(t) \|^2_{L^2(\mathcal{X})}. \label{eq:tdvp}
\end{equation}
TDVP can be interpreted as linearly approximating $u_{\text{target}}$ near the current parameter $\theta_n$ on the manifold $\mathcal{M}_{\Theta}$:
\begin{equation}
\hat{u}_{\theta_{n+1}} = \hat{u}_{\theta_n} + \nabla_\theta \hat{u}_{\theta_n} \Delta \theta + O(\|\Delta \theta\|^2). \label{eq:linear_approx}
\end{equation}
Substituting into \eqref{eq:obti} with $L^2$ loss and using \eqref{eq:target_euler} yields the least-squares problem for $\Delta\theta$:
\begin{equation}
\Delta\theta = \arg\min_{\Delta\theta} \| \Delta t \, f(\hat{u}_{\theta_n}) - \nabla_\theta\hat{u}_{\theta_n} \Delta\theta \|_{L^2}^2. \label{eq:delta_theta_ls}
\end{equation}
Dividing by $\Delta t$ and taking $\Delta t \to 0$ recovers \eqref{eq:tdvp}.
Based on this observation, TENG defines $\Delta u = u_{\text{target}} - \hat{u}_{\theta_n}$ and solves
\begin{equation}
\Delta \theta = \arg\min_{\Delta \theta} \| \Delta u - \nabla_{\theta} \hat{u}_{\theta_n} \Delta \theta \|_{L^2}, \label{eq:teng_core}
\end{equation}
with update $\theta_{n+1} = \theta_n + \Delta \theta$. This process can be iterated within each time step:
\begin{equation}
\int_{\Omega} \nabla_{\theta} \hat{u}_{\theta^{(i)}} \, dx \, \Delta \theta^{(i)} = \int_{\Omega} \Delta u^{(i)} \, dx, \label{eq:teng_iter}
\end{equation}
where superscript $(i)$ denotes sub-iteration, and $\Delta u^{(i)} = u_{\text{target}} - \hat{u}_{\theta^{(i)}}$.
%Let $\varepsilon_p(t)$ denote the $L^p$ error between the TENG solution and the exact solution. For TENG-Euler,
%\begin{equation}
%\underbrace{\|u-u_{\text{TENG}}\|}_{\varepsilon_p} \leq \underbrace{\|u-u_{\text{Euler}}\|}_{\varepsilon^{\text{EE}}} + \underbrace{\|u_{\text{Euler}}-u_{\text{TENG}}\|}_{\varepsilon^{\text{TE}}}. \label{eq:error_decomp}
%\end{equation}
%
%\begin{thm}[Chen et al. \cite{chen2024teng}]
%$\varepsilon^{\text{EE}}_p(t) = \mathcal{O}(\Delta t)$, and with $\mathcal{G} := 1 + \Delta t f$,
%\begin{equation}
%\varepsilon^{\text{TE}}_p(t) = \left\|\sum_{n=0}^{t/\Delta t-1} \mathcal{G}^n r(\cdot, t-n\Delta t)\right\|_{L^p}. \label{eq:te_error}
%\end{equation}
%\end{thm}
TENG maps physical evolution into parameter space, guiding updates via natural gradients rather than blind search. This ensures updates respect physical laws while efficiently converging to $u_{\text{target}}$.

\subsection{Relaxation Runge-Kutta Method for Structure-Preserving Targets}

The target solution $u_{\text{target}}$ is obtained through numerical discretization using classical time-stepping schemes. Whether $\hat{u}_{\theta}$ preserves physical structures depends entirely on whether $u_{\text{target}}$ retains these properties.

Consider an $s$-stage Runge-Kutta method:
\begin{align}
y_i &= u^n + \Delta t \sum_{j=1}^{s} a_{ij} f(t_n + c_j\Delta t, y_j), \quad i = 1,\ldots,s, \label{eq:rk_stages} \\
u^{n+1} &= u^n + \Delta t \sum_{i=1}^{s} b_i f(t_n + c_i\Delta t, y_i). \label{eq:rk_update}
\end{align}
Let $d^n = \sum_{i=1}^{s} b_i f_i$, so $u^{n+1} = u^n + \Delta t \, d^n$.
All Runge-Kutta methods preserve linear invariants. If coefficients satisfy
\begin{equation}
b_i a_{ij} + b_j a_{ji} = b_i b_j, \quad i,j = 1,\ldots,s, \label{eq:quadratic_condition}
\end{equation}
the method also preserves quadratic invariants. However, standard explicit Runge-Kutta methods fail to preserve quadratic or higher-order nonlinear invariants. Implicit methods capable of such preservation are computationally prohibitive in automatic differentiation frameworks due to the lack of explicit Jacobians.

To address this limitation, Ranocha \cite{ranocha2020relaxation} introduced a relaxation technique. A relaxation factor $r_n$ adaptively adjusts the step size to ensure $I(u^{n+1}_r) = I(u^n)$:
\begin{equation}
u^{n+1}_r = u^n + r \Delta t \, d^n, \label{eq:relaxed_update}
\end{equation}
where $r$ solves $F(r) = I(u^n + r \Delta t d^n) - I(u^n) = 0$. This scalar nonlinear equation can be solved efficiently via Newton's method. Under appropriate assumptions \cite{ranocha2020general}, a unique solution exists and satisfies $r = 1 + \mathcal{O}(\Delta t^{p-1})$, preserving the original accuracy order.
Integrating relaxation into TENG, we construct the target solution as
\begin{equation}
u_{\text{target}}^r = \hat{u}_{\theta_n} + r_n \Delta t \, f(\hat{u}_{\theta_n}), \label{eq:relaxed_target}
\end{equation}
ensuring $I(u_{\text{target}}^r) = I(\hat{u}_{\theta_n})$ while maintaining accuracy.

\subsection{Projection Algorithm for Parameter Updates}

Even with a structure-preserving target, the nonlinear parameterization from physical increment $\Delta u$ to parameter update $\Delta \theta$ may disrupt conservation. To ensure $I(\hat{u}_{\theta_{n+1}}) = I(\hat{u}_{\theta_n})$, we incorporate projection methods from geometric numerical integration \cite{bruna2024neural, schwerdtner2024nonlinear}.

Let $g = \nabla_{\theta} I(\theta)$ denote the gradient of the conserved quantity. Projecting the gradient onto the orthogonal complement of $g$ yields
\begin{equation}
\nabla_{\theta} \hat{u}^{\perp} = \nabla_{\theta} \hat{u} - \frac{\langle \nabla_{\theta} \hat{u}, g \rangle}{\langle g, g \rangle} g. \label{eq:projected_gradient}
\end{equation}
Substituting into the update equation gives
\begin{equation}
\nabla_{\theta} \hat{u}^{\perp} \Delta \theta = \Delta u, \label{eq:constrained_update}
\end{equation}
confining updates to the subspace orthogonal to $g$.

However, due to nonlinearity and finite network expressivity, even projected updates may drift from the conservation manifold. Following \cite{schwerdtner2024nonlinear}, we apply a post-update projection:
\begin{equation}
\theta_{I}^*(t_i) = \arg\min_{\eta \in \Theta} \frac{1}{2} \|\eta - \theta(t_i)\|^2_2 \quad \text{subject to} \quad \hat{u}(\eta,\cdot) \in \mathcal{M}_I, \label{eq:projection_optimization}
\end{equation}

Let $\mathcal{M}_I$ denote the conservation manifold. The optimal parameter set $\theta_{I}^*(t_i)$ is obtained by projecting $\theta(t_i)$ onto this manifold. Formally, we define the discrete constrained manifold in the function space as:
\begin{equation}
\hat{\mathcal{M}}_I = \left\{ \hat{u}(\eta, \cdot) \mid \eta \in \Theta, \; \hat{I}_j(\hat{u}(\eta, \cdot)) = c_j, \; j = 1, \dots, n_I \right\}, \label{eq:discrete_manifold}
\end{equation}

we solve via Lagrange multipliers. Define
\begin{equation}
\hat{m}(\eta) = [\hat{I}_1(\hat{u}(\eta,\cdot)) - \hat{I}_1(\hat{u}(\theta(0),\cdot)), \ldots, \hat{I}_{n_I}(\hat{u}(\eta,\cdot)) - \hat{I}_{n_I}(\hat{u}(\theta(0),\cdot))]^{\top}. \label{eq:constraint_function}
\end{equation}
The Lagrangian is
\begin{equation}
\mathcal{L}(\eta,\lambda) = \frac{1}{2}\|\eta - \theta_{k+1}\|_2^2 + \lambda \cdot \hat{m}(\eta). \label{eq:lagrangian}
\end{equation}
First-order optimality conditions yield
\begin{align}
\eta &= \theta_{k+1} + \hat{m}'(\eta)^{\top} \lambda, \\
0 &= \hat{m}(\eta). \label{eq:kkt}
\end{align}
The Lagrange multiplier is updated iteratively:
\begin{align}
\Delta\lambda^{(i)} &= -(\hat{m}'(\theta_{k+1})\hat{m}'(\theta_{k+1})^{\top})^{-1} \hat{m}(\theta_{k+1} + \hat{m}'(\theta_{k+1})^{\top} \lambda^{(i)}), \\
\lambda^{(i+1)} &= \lambda^{(i)} + \Delta\lambda^{(i)}. \label{eq:lagrange_update}
\end{align}
It should be noted that the execution of the projection algorithm, when implementing the Lagrange multiplier method, is susceptible to the singularity of the Jacobian matrix $\nabla_{\theta} \hat{u}$

\subsection{The Complete RP-TENG Algorithm}

Combining both techniques, we obtain the Relaxation-Projection TENG (RP-TENG) method, summarized in Algorithm \ref{alg:rpteng}.

\begin{algorithm}[H]
\caption{Relaxation-Projection TENG (RP-TENG) with Euler base scheme}
\label{alg:rpteng}
\begin{algorithmic}[1]
\Require Time step $\Delta t$, final time $T_{\text{End}}$, network architecture $\hat{u}_{\theta}$, initial condition $u_0$, right-hand side $f$
\State $\theta_0 \gets \arg\min_{\eta} \|\hat{u}_{\eta} - u_0\|^2$ \Comment{Initial parameter fitting}
\State $t \gets 0$, $\theta_{\text{cur}} \gets \theta_0$, $\hat{u}_{\text{cur}} \gets \hat{u}_{\theta_0}$
\While{$t < T_{\text{End}}$}
    \State Compute $f_{\text{cur}} = f(\hat{u}_{\text{cur}})$ via automatic differentiation
    \State $u_{\text{temp}} \gets \hat{u}_{\text{cur}} + \Delta t \, f_{\text{cur}}$
    \State Find $\gamma$ solving $I(u_{\text{temp}}(\gamma)) = I(\hat{u}_{\text{cur}})$ \Comment{Relaxation step}
    \State $u_{\text{target}} \gets \hat{u}_{\text{cur}} + \gamma \Delta t \, f_{\text{cur}}$
    \State $\theta_{\text{temp}} \gets \text{TENG\_update}(\theta_{\text{cur}}, u_{\text{target}})$ \Comment{Solve \eqref{eq:teng_core}}
    \State $\theta_{\text{new}} \gets \text{project\_onto\_manifold}(\theta_{\text{temp}})$ \Comment{Solve \eqref{eq:projection_optimization}}
    \State $\theta_{\text{cur}} \gets \theta_{\text{new}}$, $\hat{u}_{\text{cur}} \gets \hat{u}_{\theta_{\text{new}}}$
    \State $t \gets t + \gamma \Delta t$
\EndWhile
\State \Return $\theta_{\text{cur}}$
\end{algorithmic}
\end{algorithm}

\section{Numerical Experiments}

This section presents numerical experiments demonstrating the effectiveness of RP-TENG on three benchmark problems. All experiments are implemented in Python using JAX \cite{jax2018github} for automatic differentiation with double-precision arithmetic. Reference solutions are obtained via spectral methods. For the solution trajectory at times $\{t_k\}_{k=0}^n$, we compute relative errors:
\begin{equation}
E_r(t_k) = \frac{\sum_{i=1}^{n_E} \|\hat{u}(\theta_k, x_i^{\text{test}}) - u^{\text{ref}}(t_k, x_i^{\text{test}})\|}{\sum_{i=1}^{n_E} \|u^{\text{ref}}(t_k, x_i^{\text{test}})\|}, \label{eq:relative_error}
\end{equation}
where $\{x_i^{\text{test}}\}_{i=1}^{n_E}$ are equally spaced test points.

\subsection{Inviscid Burgers Equation: Mass Conservation}

We first consider the inviscid Burgers equation on $[-1,1)$ with periodic boundary conditions:
\begin{equation}
\partial_t u + u \partial_x u = 0, \quad u(x,0) = 1 + 0.3 \exp(-B^2 x^2), \label{eq:burgers}
\end{equation}
where $B$ controls the Gaussian width. The equation conserves mass: $I_{\text{mass}} = \int_{-1}^1 u \, dx$.

The solution is parameterized by a fully-connected neural network consisting of five hidden layers, each with 10 neurons and $\tanh$ activation, totaling $n_{\bm{\theta}} = 371$ trainable parameters. Periodic boundary conditions are strictly enforced through input-layer encoding \cite{berman2023randomized}. To mitigate numerical instabilities arising from the potential singularity of $\nabla_{\bm{\theta}}\hat{u}$, a randomized update strategy \cite{berman2023randomized} is employed. We utilize $n_S = 1000$ equidistant sample points for the least-squares problem \eqref{eq:teng_core}, while $n_M = 1000$ quadrature points are used for the evaluation of energy errors.

Define mass error as

\begin{equation}
E_{\text{mass}} = |I(\hat{u}_{\theta(0)}) - I(\hat{u}_{\theta_n})|
\end{equation}

\begin{figure}[htbp]
    \centering
    \includegraphics[width=0.8\textwidth]{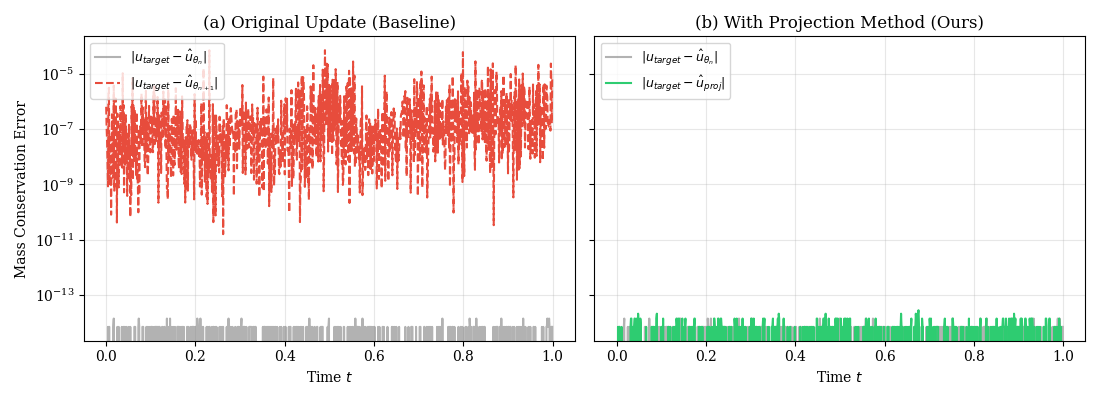} % 替换为你的文件名
    \caption{Evolution of the mass conservation error over time. The projection method (blue solid line) maintains the error near machine precision ($\sim 10^{-14}$), significantly outperforming the original method (orange dashed line) which exhibits fluctuations between $10^{-9}$ and $10^{-5}$.}
    \label{fig:mass_conservation}
\end{figure}

For the inviscid Burgers' equation, mass is a quintessential linear conserved quantity. Although classical Runge-Kutta methods are theoretically known to preserve linear conservation laws, a significant conservation error is observed during the neural numerical solving process between the updated network solution $\hat{u}_{\theta_{n+1}}$ and the target solution $u_{\text{target}}$. As illustrated in Figure 1(a), this error accumulates over time, leading to a severe deviation of the numerical solution from the initial mass.
In contrast, the target solution $u_{\text{target}}$ strictly constrains the mass error relative to the network solution at the previous time step $\hat{u}_{\theta_n}$ within the order of $10^{-13}$. Figure 1(b) further validates the necessity of the projection algorithm within the TENG framework: the solution $\hat{u}_{\theta_{\text{proj}}}$, obtained through the projection operator, perfectly preserves the mass characteristics of $\hat{u}_{\theta_n}$. This provides strong evidence for the effectiveness of our proposed projection algorithm in maintaining long-term physical conservation.

% 在导言区（\begin{document}之前）加入这个设置
% 正文排版

\begin{figure}[htbp]
    \centering
    \includegraphics[width=1.0\textwidth]{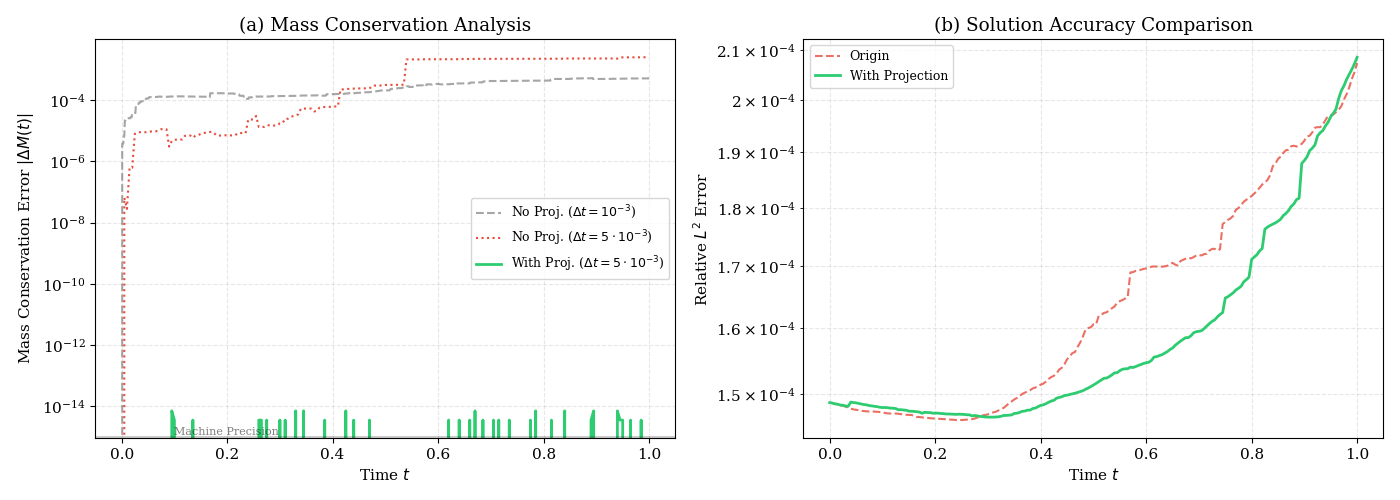}
    \caption{Numerical performance of the TENG framework on the inviscid Burgers' equation.
    \textbf{(a)} Mass conservation error relative to different time step sizes and the projection method.
    \textbf{(b)} Temporal evolution of the relative $L^2$ error, highlighting the accuracy gains from the projection-based physical constraint.}
    \label{fig:22}
\end{figure}

We evaluated the mass conservation error relative to the initial state using $n_M = 200$ quadrature points for time steps $\Delta t = 5 \times 10^{-3}$ and $10^{-3}$. The results Fig 2 indicate that the vanilla TENG framework fails to preserve the linear invariant, even with a refined temporal resolution. However, by incorporating the projection algorithm at $\Delta t = 5 \times 10^{-3}$, the mass error is effectively suppressed. The numerical evidence clearly demonstrates that the projection step ensures physical consistency and prevents the long-term drift characteristic of unconstrained neural solvers.

Figure \ref{fig:22} provides a comprehensive evaluation of the TENG framework. Subfigure \textbf{(a)} demonstrates that while simply reducing the time step $\Delta t$ is insufficient to rectify the mass drift inherent in unconstrained neural solvers, the proposed projection method effectively constrains the error to near-machine precision. This geometric enforcement directly translates into sustained predictive accuracy, as evidenced by the relative $L^2$ error shown in \textbf{(b)}. These results suggest that maintaining physical invariants is a prerequisite for achieving high-fidelity, long-term neural simulations.

\subsection{Korteweg-de Vries Equation: Quadratic Invariant Preservation}

The Korteweg-de Vries (KdV) equation,
\begin{equation}
\partial_t u + \partial_x\left(\frac{u^2}{2}\right) + \partial_x^3 u = 0, \label{eq:kdv}
\end{equation}
admits single-soliton solutions
\begin{equation}
u(x,t) = A \cosh^{-2}\left(\frac{\sqrt{3A}(x-ct-\mu)}{6}\right), \quad c = \frac{A}{3}, \label{eq:kdv_soliton}
\end{equation}
with $A=2$ in our experiments. Conserved quantities include mass $I_{\text{mass}} = \int u \, dx$ and energy $I_{\text{energy}} = \int u^2 \, dx$ (quadratic invariant).

The governing equation is solved on the spatial domain $\Omega = (0, 80)$ with periodic boundary conditions enforced via a periodic embedding layer at the input stage. The core architecture is a fully connected MLP with three hidden layers (10 neurons each), employing sinusoidal activation $\sigma(x) = \sin(x)$ to capture the high-order derivatives and oscillatory nature of the KdV solutions. The model comprises 371 trainable parameters and utilizes $n_S = 1000$ integration nodes. Temporal integration is performed using an explicit RK4 scheme with a fixed step $\Delta t = 5 \times 10^{-3}$. To assess physical fidelity, the energy conservation error is quantified as $E_{\text{err}}(t_n) = |I(\hat{u}_{\theta_0}) - I(\hat{u}_{\theta_n})|$.

\begin{figure}[htbp]
    \centering
    \includegraphics[width=1.0\textwidth]{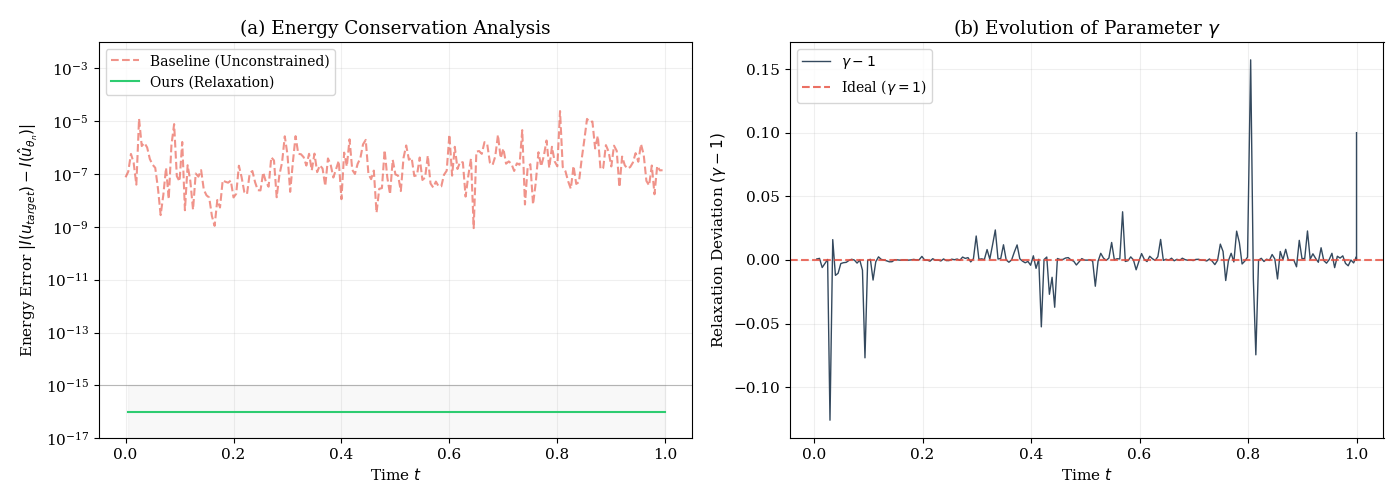}
\caption{Numerical performance of the TENG framework on the KdV equation.
\textbf{(a) Necessity of the relaxation method:} The RP-TENG framework suppresses the quadratic invariant error $|\mathcal{H}(u_{\mathrm{target}}) - \mathcal{H}(\hat{u}_{\theta_{n}})|$ to machine epsilon ($\sim 10^{-15}$), effectively eliminating the numerical drift of the baseline.
\textbf{(b) Influence on dynamical trajectory:} The stability of $\gamma \approx 1$ confirms that the energy constraint is enforced with minimal intrusiveness to the neural dynamical trajectory.}
    \label{fig:combined}
    \label{fig:combined}
\end{figure}

Previous tests on the inviscid Burgers' equation failed to fully showcase the advantages of the relaxation method, as the mass is a linear invariant naturally preserved by standard Runge-Kutta schemes. To provide a more rigorous assessment, we investigate the Korteweg-de Vries (KdV) equatio. Its energy, a quadratic invariant, is notably difficult for conventional explicit integrators to maintain. This scenario offers a more demanding benchmark to demonstrate the energy-conserving capabilities of the RP-TENG framework.

Consequently, we present a comparative analysis between the vanilla and the relaxed method. As illustrated in  Figure 3(a), the relaxation algorithm effectively constrains the numerical solution to the energy manifold, ensuring the physical fidelity of $\hat{u}_{\theta_{n}}$. Furthermore,Figure 3(b) depicts the range of the relaxation parameter $\gamma$. The fact that $\gamma$ remains concentrated near unity indicates that the proposed algorithm imposes physical constraints without significantly perturbing the intrinsic temporal dynamics of the system.

\begin{figure}[htbp]
    \centering
    \includegraphics[width=1.0\textwidth]{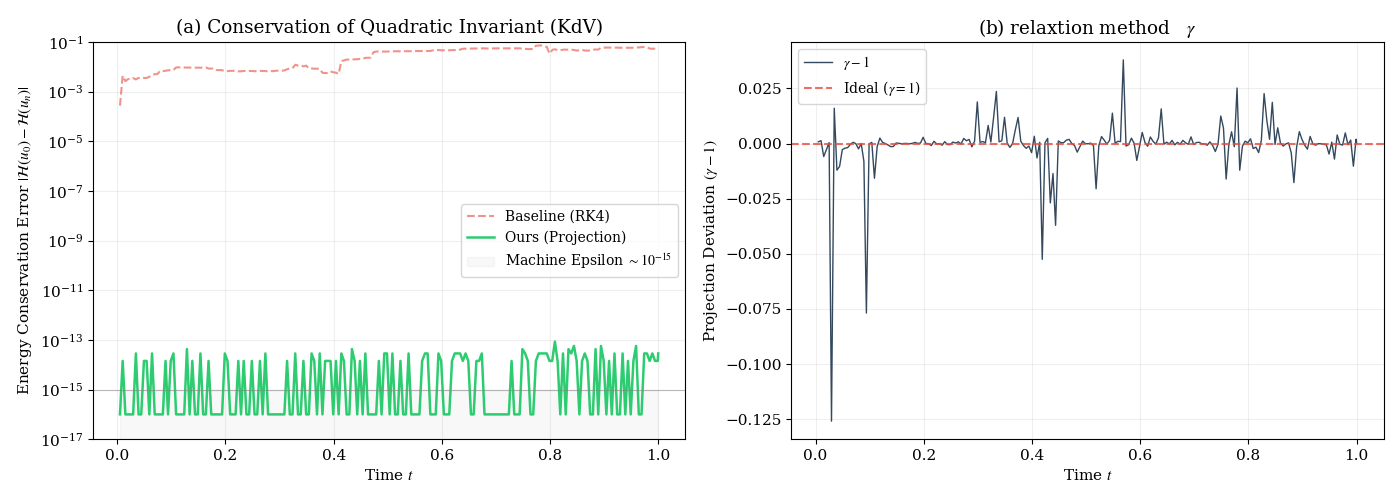}
\caption{Numerical performance of RP-TENG on the KdV equation.
    \textbf{(a) Energy conservation:} RP-TENG maintains the Hamiltonian error $|\mathcal{H}(u_0) - \mathcal{H}(\hat{u}_{\theta_{n}})|$ near machine epsilon, eliminating the drift of vanilla TENG.
    \textbf{(b) Parameter stability:} The temporal consistency of $\gamma \approx 1$ confirms the non-intrusiveness of the projection step within the dual-constraint framework.}
    \label{fig:3}
\end{figure}

\begin{figure}[htbp]
    \centering
    \begin{minipage}[t]{0.47\textwidth}
        \centering
        \includegraphics[width=\textwidth]{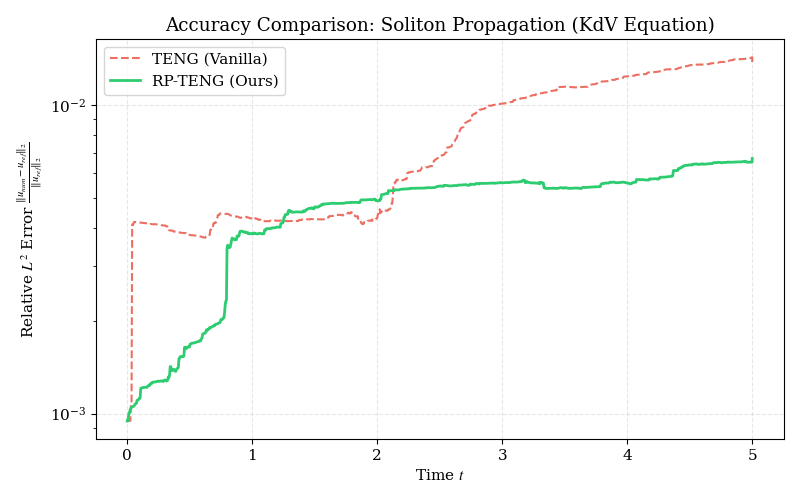}
        \centerline{\textbf{(a) Solution Accuracy }}
    \end{minipage}
    \hfill
    \begin{minipage}[t]{0.48\textwidth}
        \centering
        \includegraphics[width=\textwidth]{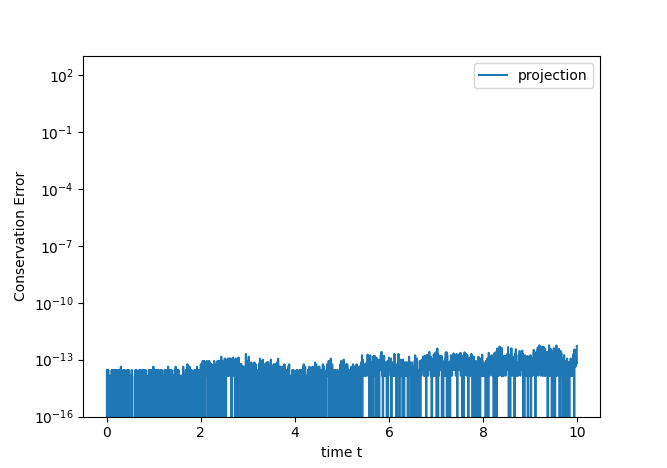}
        \centerline{\textbf{(b)   Energy Conservation}}
    \end{minipage}
    \caption{\textbf{Comprehensive performance on the KdV equation soliton simulation.} (a) Temporal evolution of the relative $L^2$ error against the analytical solution.(b) Long-term energy-preserving capability of the RP-TENG framework for $t \in [0, 10]$ }
    \label{fig:5}
\end{figure}

 Figure 4(a) highlights the superior energy-preserving properties of the RP-TENG framework. For the KdV equation, the quadratic energy invariant is highly susceptible to numerical dissipation; however, the integrated Relaxation-Projection mechanism effectively suppresses Hamiltonian drift, ensuring physical consistency during long-term soliton propagation. To assess the impact of these constraints, Figure 4(b) examines the evolution of the relaxation scaling factor $\gamma$. The marginal fluctuations around unity (deviations $\sim 10^{-12}$) validate the non-intrusiveness of the adaptive correction, confirming that the framework enforces physical invariants with minimal perturbation to the intrinsic dynamical trajectory.

As shown in  Figure \ref{fig:5}, RP-TENG demonstrates superior numerical robustness over the vanilla TENG. Specifically, the relative $L_2$ error growth is effectively mitigated Figure 5(a), as the manifold-anchoring mechanism prevents the phase drift and artificial dissipation common in unconstrained models. Furthermore, while the vanilla TENG exhibits significant energy divergence, RP-TENG restricts Hamiltonian drift to machine precision ($\sim 10^{-15}$) throughout the simulation Figure 5(b). These results confirm that enforcing physical invariants is essential for maintaining the long-term fidelity of neural PDE solvers in nonlinear regimes.

\subsection{Acoustic Wave Equation: Hamiltonian Preservation}

Consider the first-order acoustic wave system:
\begin{align}
\frac{\partial \rho}{\partial t} &= -\frac{\partial v}{\partial x}, \\
\frac{\partial v}{\partial t} &= -\frac{\partial \rho}{\partial x}, \label{eq:acoustic}
\end{align}
on $[-1,1)$ with periodic boundaries,initial conditions $\rho(x,0) = e^{-9x^2}$, $v(x,0)=0$. The Hamiltonian (total energy) is conserved:
\begin{equation}
H_{\text{wave}}(\rho,v) = \frac{1}{2} \int_{\mathcal{X}} (\rho(x)^2 + v(x)^2) \, dx. \label{eq:hamiltonian}
\end{equation}

Numerical implementation  for the acoustic wave equation incorporates a periodic embedding layer  at the input stage to naturally satisfy the boundary constraints. The underlying architecture is a fully connected multi-layer perceptron (MLP) consisting of four hidden layers, each containing 20 neurons, with a two-dimensional output. The $\tanh$ activation function is employed throughout the hidden layers to ensure smooth representation of the wave fields. Performance evaluation of the PR-TENG framework relies on a high-fidelity reference solution synthesized via the Fourier spectral method for spatial discretization (over 256 grid points) and the RK45 scheme for temporal integration ($\Delta t = 10^{-3}$). Consequently, the prediction accuracy is quantified by the relative $L^2$ error as defined in (\ref{eq:relative_error}).

\begin{figure}[htbp]
    \centering
    \includegraphics[width=1.0\textwidth]{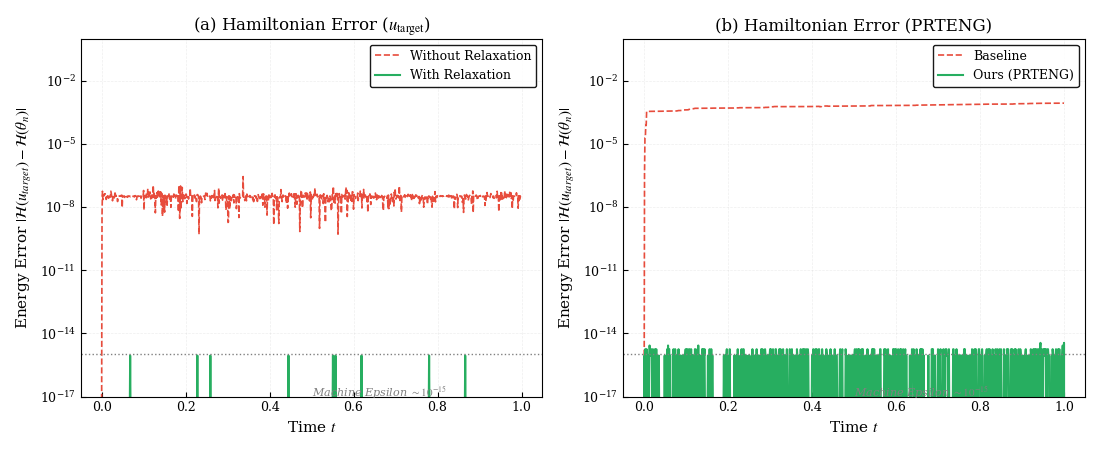}
    \caption{Numerical performance of the PR-TENG framework on  the Acoustic Wave Equation
\textbf{(a) } Comparison of the Hamiltonian error for $u_{\mathrm{target}}$ with and without the relaxation scheme, illustrating the mitigation of energy drift.
    \textbf{(b)} Comparative performance of the PR-TENG and the original TENG algorithms.}
    \label{fig:6}
\end{figure}

\begin{figure}[htbp]
    \centering
    \includegraphics[width=0.6\textwidth]{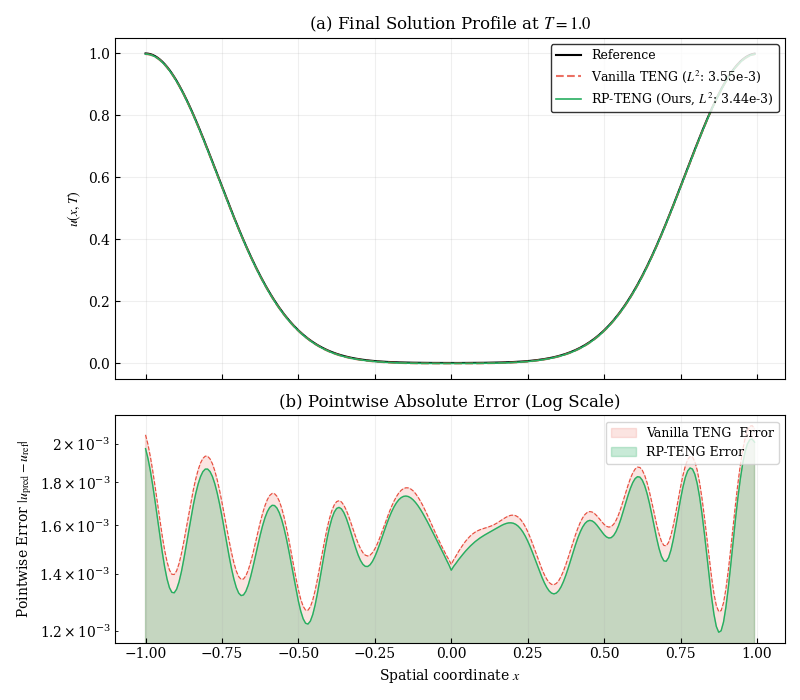}
    \caption{Relative error of RP-TENG on the acoustic wave equation}
    \label{fig:7}
\end{figure}

\begin{figure}[htbp]
    \centering
    \includegraphics[width=0.6\textwidth]{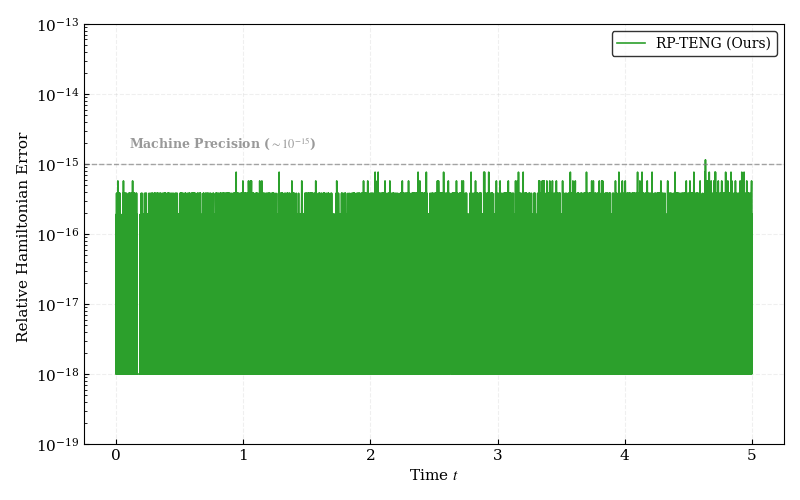}
    \caption{ RP-TENG long time implement effects }
    \label{fig:8}
\end{figure}

\ref{fig:6} evaluates the integration of the relaxation method within the TENG framework. The execution of the relaxation scheme enables the acquisition of the target solution $u_{\mathrm{target}}$ while successfully preserving the system's Hamiltonian energy, thereby demonstrating the necessity of the proposed relaxation approach. Furthermore, tests conducted over the temporal interval $t \in [0, 1]$ reveal that the RP-TENG algorithm ensures the strict conservation of the Hamiltonian invariant, confirming its robustness in maintaining the underlying physical properties of the system.

Figure \ref{fig:7} illustrates the relative errors of the vanilla TENG and the proposed RP-TENG framework at $T=1$ with a time step of $\Delta t = 10^{-3}$. For the acoustic wave equation, the RP-TENG achieves a relative error of $3.44 \times 10^{-3}$, which is a noticeable improvement over the $3.55 \times 10^{-3}$ produced by the original TENG. This result indicates that the RP-TENG framework not only maintains the Hamiltonian invariant through the relaxation and projection mechanisms but also enhances the predictive accuracy compared to the baseline TENG scheme.

Figure 8 evaluates the energy-preserving performance of the RP-TENG algorithm applied to the acoustic wave equation. With a temporal discretization of $\Delta t = 10^{-3}$, the numerical integration is conducted over an extended temporal domain $t \in [0, 5]$

\section{Conclusion}

Traditional numerical analysis emphasizes the preservation of the physical properties inherent in the original system. As an emerging paradigm for solving partial differential equations (PDEs), neural networks are naturally expected to maintain these intrinsic energy attributes. Given that neural network outputs are dictated by their parameter configurations, and that physical properties are frequently time-dependent, conventional global training frameworks often prove inadequate. In contrast, a local training framework  is better suited to capturing the temporal dynamics of the system, offering a distinct advantage in structure-preserving tasks.

In this work, we enhance the \textit{Time-evolution Neural Gradient} (TENG) method by integrating two synergistic structure-preserving mechanisms, resulting in the RP-TENG framework:
 The relaxation method  that adaptively re-scales the update to ensure the target solution $u_{\text{target}}$ strictly adheres to the prescribed invariants;
  and the projection algorithm  that maps the evolved network parameters back onto the underlying conservation manifold.
The resulting RP-TENG framework provides the following key advantages:

\begin{itemize}
    \item \textbf{Structure-preserving targets:} The relaxation scheme ensures that $u_{\text{target}}$ inherits conservation properties from the preceding time step without compromising the formal order of accuracy.
    \item \textbf{Constrained manifold updates:} By combining tangent-space projection with post-update correction, the framework keeps the trajectory of the network parameters strictly on the invariant manifold.
    \item \textbf{Enhanced long-term robustness:} Numerical experiments on the Burgers, KdV, and acoustic wave equations demonstrate that RP-TENG significantly mitigates non-physical energy drift while maintaining high predictive precision over extended temporal domains.
\end{itemize}

Future research will focus on extending this framework to dissipative systems (e.g., phase-field models) and exploring its applicability to stochastic differential equations governed by latent invariants.
\bibliographystyle{plain}
\bibliography{ref}

\end{document}